\documentclass[12pt]{amsart}

\usepackage{fullpage}

\usepackage{amssymb}
\usepackage{amsthm}
\usepackage{amsmath}
\usepackage{amsrefs}

\newtheorem{theorem}{Theorem}

\newtheorem{lemma}{Lemma}
\newtheorem{proposition}{Proposition}

\theoremstyle{remark}

\newcommand{\N}{\mathbb{N}}
\newcommand{\R}{\mathbb{R}}
\newcommand{\ped}{\mathrm{Ped}}
\newcommand{\Ideal}{\mathrm{Ideal}}
\newcommand{\precp}{\preceq_{CP}}
\newcommand{\simcp}{\sim_{CP}}

\begin{document}

\title[On the comparison of positive elements of a C*-algebra...]
{On the comparison of positive elements of a C*-algebra by lower semicontinuous traces}

\author{Leonel Robert}

\address{Leonel Robert, Fields Institute, Toronto, Canada, M5T 3J1}

\email{lrobert@fields.utoronto.ca}


\begin{abstract}
It is shown in this paper that two positive elements of a C*-algebra agree on all
lower semicontinuous traces if and only if they are equivalent in the sense of Cuntz
and Pedersen. A similar result is also obtained in the more general case where the two elements
are comparable by their values on the lower semicontinuous traces.
This result is used to give a characterization
of the functions on the cone of lower semicontinuous traces of a stable C*-algebra that arise from 
positive elements of the algebra.   
\end{abstract}
\maketitle

\section{Introduction}
In \cite{cuntz-pedersen}, Cuntz and Pedersen considered the problem of comparing 
positive elements
of a C*-algebra by their values on the lower semicontinuous traces on the
algebra. They defined an equivalence relation on the positive elements---the Cuntz-Pedersen
relation---and showed, among other results, that if the C*-algebra
is simple, then two positive elements are Cuntz-Pedersen equivalent 
if and only if they agree on all the lower semicontinuous traces. The question whether this was true for an
arbitrary C*-algebra was left unsettled in their paper, and is answered affirmatively
in Theorem \ref{maincomparison} below.

Let $A$ be  C*-algebra. Recall that a trace on $A$ is a map $\tau\colon A^+\to
[0,\infty]$ that is additive, homogeneous, and satisfies the 
identity $\tau(xx^*)=\tau(x^*x)$. We will be mostly  interested in the lower
semicontinuous traces on $A$, the set of which we shall denote by $T(A)$.

For $a,b\in A^+$ let us say that $a$ is
Cuntz-Pedersen equivalent to $b$, denoted by $a\sim b$, if $a=\sum_{i=1}^\infty x_{i}x_i^*$
and $b=\sum_{i=1}^\infty x_i^*x_i$
for some $x_i\in A$. Let us say that $a$ is Cuntz-Pedersen smaller than $b$,
and denote this by $a\preceq b$, if $a\sim a'\leq b$ for some $a'\in A^+$. It was shown in \cite{cuntz-pedersen} 
that the relations $\sim$ and $\preceq$ are transitive (so $\sim$ is an equivalence relation).

\begin{theorem}\label{maincomparison}
Let $a$ and $b$ be positive elements of $A$. The following propositions are
true.

(i) $\tau(a)\leq \tau(b)$ for all $\tau\in T(A)$ if and only if for every
$\epsilon>0$ there is $\delta>0$ such that
$(a-\epsilon)_+\preceq (b-\delta)_+$.

(ii) $\tau(a)=\tau(b)$ for all $\tau\in T(A)$ if and only if $a\sim b$.
\end{theorem}

\emph{Remark.} We shall see in Section 3 an example where $\tau(a)\leq \tau(b)$ for all $\tau\in T(A)$,
but it is not true that $a\preceq b$.

\section{Proof of Theorem \ref{maincomparison}}
The proof of Theorem \ref{maincomparison} is preceded
by a number of preliminary definitions and results.

Let $a,b\in A^+$. Let us write $a\precp b$ if for every $\epsilon>0$ there
is $\delta>0$
such that $(a-\epsilon)_+\preceq (b-\delta)_+$. Since $\preceq$ is
transitive, the relation
$\precp$ is clearly transitive as well. Let us write $a\simcp b$ if 
$a\precp b$ and $b\precp a$. This defines an equivalence
relation in $A^+$. (It will be shown in the proof of
Theorem \ref{maincomparison} (ii) that $\simcp$ is the same as $\sim$.) 

\begin{proposition}
Let $a,b,c,d\in A^+$. The following propositions are true.

(i) If $a\precp b$ and $c\precp d$ then $a+c\precp b+d$.

(ii) If $a\precp b$ then $\alpha a\precp \alpha b$ for all $\alpha\in \R^+$.

(iii) If $a\preceq b$ then $a\precp b$.
\end{proposition}
\begin{proof}
It was shown in \cite[Proposition 2.3]{elliott-robert-santiago} that for all $a,b\in A^+$ and $\epsilon>0$
there is $\delta>0$ such that
\begin{align}
(a-\epsilon)_+ +(b-\epsilon)_+ &\preceq (a+b-\delta)_+,\label{CP1}\\
(a+b-\epsilon)_+ &\preceq (a-\delta)_+ + (b-\delta)_+\label{CP2}.
\end{align}
These inequalities imply (i).

(ii) This is clear.

(iii) Suppose that $a=\sum_{i=1}^\infty x_ix_i^*$ and $\sum_{i=1}^\infty x_i^*x_i\leq b$.
Let $\epsilon>0$. By the lemma of Kirchberg and R\o rdam \cite[Lemma 2.2]{kirchberg-rordam}
(see also the remark after \cite[Lemma 2.2]{elliott-robert-santiago}), there
are $n\in \N$ and $\epsilon_1>0$ such that $(a-\epsilon)_+\preceq (\sum_{i=1}^n x_ix_i^*-\epsilon_1)_+$.
We have
\begin{align*}
(a-\epsilon)_+\preceq (\sum_{i=1}^n x_ix_i^*-\epsilon_1)_+ \preceq \sum_{i=1}^n(x_ix_i^*-\epsilon_2)_+
&\sim \sum_{i=1}^n(x_i^*x_i-\epsilon_2)_+\\
&\preceq (\sum_{i=1}^n x_ix_i^*-\epsilon_3)_+\preceq (b-\epsilon_4)_+. 
\end{align*}
In the above chain of inequalities we have applied \eqref{CP1}, \eqref{CP2} and that 
$(xx^*-\epsilon)_+\sim (x^*x-\epsilon)_+$ for all $\epsilon>0$ (see \cite[Proposition 2.3]{elliott-robert-santiago}).
\end{proof}

Let us denote by $A_{CP}$ the quotient $A^+/\!\simcp$. We consider $A_{CP}$ ordered
by the order
$\langle a\rangle\leq \langle b\rangle$ if $a\precp b$, where $\langle a\rangle$ and $\langle b\rangle$ denote the equivalence
classes of the positive elements $a$ and $b$.
We also endow $A_{CP}$ with the addition operation $\langle a\rangle+\langle b\rangle:=\langle a+b\rangle$. 
The order of $A_{CP}$ is compatible with the addition operation, i.e.,
$\langle a\rangle\leq \langle a\rangle+\langle b\rangle$. 
Thus, $A_{CP}$ is an ordered semigroup with $0$. 

In order to prove Theorem \ref{maincomparison} (i) we will apply the
following proposition to the ordered semigroup $A_{CP}$.

\begin{proposition}\label{strictcomp} (\cite[Proposition 3.2]{rordam}) 
Let $S$ be an ordered semigroup with 0 and with the property that if $(k+1)x\leq
ky$ for some $x,y\in S$ and $k\in \N$,
then $x\leq y$ (i.e., $S$ is almost unperforated). The following implication holds in $S$:

If $x\leq My$ for some $M\in \N$, and $\lambda(x)<\lambda(y)$ for every
$\lambda\colon S\to [0,\infty]$ that is additive, order-preserving,
and satisfies $\lambda(y)=1$, then $x\leq y$.
\end{proposition}

Notice that the semigroup $A_{CP}$ satisfies the hypotheses of the preceding
proposition. In fact, in $A_{CP}$ we have that $k\langle x\rangle\leq k\langle y\rangle$
implies $\langle x\rangle\leq \langle y\rangle$, because we can multiply the
elements of $A_{CP}$ by positive real scalars.

Notice also that for every additive map 
$\lambda\colon A_{CP}\to [0,\infty]$, the map on $A^+$ defined
by $\tau_\lambda(a):=\lambda(\langle a\rangle)$ is a trace, because it is 
additive and satisfies the trace identity (homogeneity holds automatically
for any additive map with values in $[0,\infty]$). This trace may not be lower semicontinuous. 
To Theorem \ref{maincomparison} (i) we will then need the following lemma.

\begin{lemma}\label{regularization}(\cite[Lemma 3.1]{elliott-robert-santiago})
Let $\tau\colon A^+\to [0,\infty]$ be  a trace on the C*-algebra $A$. Then
$\tilde\tau(a)=\sup_{\epsilon>0} \tau((a-\epsilon)_+)$
is a lower semicontinuous trace.
\end{lemma}

\begin{proof}[Proof of Theorem \ref{maincomparison} (i)]
It is clear that if $a\precp b$ then $\tau(a)\leq \tau(b)$ for every
$\tau\in T(A)$.

Let us assume that $\tau(a)\leq \tau(b)$ for all $\tau\in T(A)$. For
every closed two-sided ideal $I$ of $A$, the map defined by $\tau_I(x)=0$ if $x\in
I^+$ and $\tau_I(x)=\infty$ otherwise,
is a  lower semicontinuous trace on $A$. Since $\tau_I(a)\leq \tau_I(b)$ for any such trace
it follows that $\Ideal(a)\subseteq \Ideal(b)$. Let $\epsilon>0$. 
We have  $(a-\epsilon)_+\in \ped(\Ideal(b))^+$; hence 
$(a-\epsilon)_+=\sum_{i=1}^m{y_iby_i^*}$ for some $y_i\in A$. This implies that $(a-\epsilon)_+\preceq Mb$
for some $M\in \N$, and so $\langle(a-\epsilon)_+\rangle\leq M\langle b\rangle$. Let us show that
we also have $\lambda(\langle(a-\epsilon)_+\rangle) < \lambda (\langle b\rangle)$ for any
additive, order-preserving map $\lambda\colon A_{CP}\to [0,\infty]$ such that 
$\lambda(\langle b\rangle)=1$.
By Proposition \ref{strictcomp}, this will imply that 
$\langle(a-\epsilon)_+\rangle\leq \langle b\rangle $, and since
$\epsilon$ is arbitrary, we will have $\langle a\rangle\leq \langle b\rangle$, as desired.

Let $\lambda\colon A_{CP}\to [0,\infty]$ be additive, order-preserving, and
such that $\lambda(\langle b\rangle)=1$. Let $\tau_\lambda$ be the trace associated to
$\lambda$ and $\tilde\tau_{\lambda}$ its lower semicontinuous regularization defined
as in Lemma \ref{regularization}. Notice that for every $\delta>0$ and $c\in A^+$ we have
$\tau_\lambda((c-\delta)_+)\leq \tilde\tau_\lambda(c)\leq \tau_\lambda(c)$.

Case 1. Assume that $\tilde\tau_\lambda(a)\neq 0$. Then
\[
\lambda(\langle(a-\epsilon)_+\rangle)=\tau_\lambda((a-\epsilon)_+)\leq
\tilde\tau_\lambda((a-\epsilon/2)_+)<\tilde \tau_\lambda(a)\leq \tilde \tau_\lambda(b)\leq
\tau_\lambda(b)=\lambda(\langle b\rangle).
\]

Case 2. Assume that $\tilde\tau_\lambda(a)= 0$. Then
\[
\lambda(\langle(a-\epsilon)_+\rangle)=\tau_\lambda((a-\epsilon)_+) = 0<1=\lambda(\langle b\rangle).
\qedhere
\]
\end{proof}

We now turn to the proof of Theorem \ref{maincomparison} (ii). 

\begin{lemma}\label{cancellation}
Let $a,b,c\in A^+$. The following propositions are true.

(i) If $a+c\preceq b+c$ and $c\preceq Mb$ for some $M\in \R^+$, then $a\preceq (1+\delta)b$
for all $\delta>0$.

(ii) If $a+c\sim b+c$ and $c\preceq Ma,Mb$ for some $M\in \R^+$, then $a\sim b$.
\end{lemma}
\begin{proof}
(i) Adding $a$ on both sides of  $a+c\preceq b+c$ we obtain that $2a+c\preceq 2b+c$, and by
induction $ka+c\preceq kb+c$ for all $k$. Dividing the last inequality by $k$ we see that we may assume that
$c\preceq b$. Changing $c$ to $c'$ such that $c\sim c'\leq b$ we may further assume that
$c\leq b$. We have $a+b=a+c+(b-c)\preceq 2b$. By applying induction
we obtain that $ka+b\preceq (k+1)b$ for all $k$. This proves (i).

(ii) As in the proof of (i) we may assume that $c\leq b$. We have 
$a+b\sim a+c+b-c\sim 2b$. In the same way we obtain that $a+b\sim 2a$. Therefore, $a\sim b$.
\end{proof}

\begin{proof}[Proof of Theorem \ref{maincomparison} (ii)]
By Theorem \ref{maincomparison} (i) we know that $a\precp b$ and $b\precp
a$. Thus, there are strictly decreasing sequences of positive numbers
$\epsilon_n$ and $\mu_n$
such that $(a-\epsilon_n)_+\preceq (b-\mu_n)_+\preceq (a-\epsilon_{n+1})_+$
for all $n\in \N$, and $\epsilon_n,\mu_n\to 0$.
Let us write $a_n=(a-\epsilon_n)_+$ and $b_n=(b-\mu_n)$. The proof now follows 
the same arguments as in \cite[Lemma 7.4]{cuntz-pedersen}, although with some modifications.

Since $a_n\preceq
b_n$ there is $a_n'$ such that $a_n\sim a_n'\leq b_n$.
We have $(b_n-a_n')+a_n'\preceq (a_{n+1}-a_n) + a_n$. We also have $a_n\leq
M(a_{n+1}-a_n)$. (This is verified by referring back to
the definition of $a_n$ in terms of $a$.) By Lemma \ref{cancellation} (i) we have that $b_n-a_n'\preceq
(1+\delta)(a_{n+1}-a_n)$ for all $\delta>0$. 
For $\delta$ small enough we have that $b_n-a_n'\preceq
(1+\delta)(a_{n+1}-a_n)\leq a_{n+2}-a_n$.
Let $z$ be such that $b_n-a_n'\sim z\leq a_{n+2}-a_n$ and let $b_n'=a_n+z$.
The sequence $(b_n')$ satisfies that $b_n'\sim b_n$ and $a_n\leq b_n'\leq a_{n+2}$.

We have $(b_{n+3}'-b_n')+b_n'\sim (b_{n+3}-b_n)+b_n$. Let us see that we can
apply Lemma \ref{cancellation} (ii) to conclude that
$(b_{n+3}'-b_n')\sim  (b_{n+3}-b_n)$. On one hand, $b_n\leq M(b_{n+3}-b_n)$ for some
$M>0$,
by the expression of these elements
in terms of the functional calculus of $b$. On the other hand, $b_n'\leq
a_{n+2}\leq M'(a_{n+3}-a_{n+2})\leq M'(b_{n+3}'-b_n')$ for some $M'>0$.
Hence, Lemma \ref{cancellation} (ii) can be applied. We have
\[
a=b_0'+\sum_{k=1}^\infty (b_{3(k+1)}'-b_{3k}')\sim b_0+\sum_{k=1}^\infty
(b_{3(k+1)}-b_{3k})=b.
\qedhere
\]
\end{proof}

\section{Further remarks}
Let us first show that Theorem \ref{maincomparison} (i) cannot be strengthened to obtain 
$\tau(a)\leq \tau(b)$ for all $\tau$ implies $a\preceq b$. 

\emph{Example.} Let $A$ be a simple unital C*-algebra with exactly 
two extreme tracial states $\tau_1$ and $\tau_2$. Such an algebra may be found,
say, among the AF C*-algebras (e.g., \cite[Example 6.10]{cuntz-pedersen}). Let $b\in A^+$ be such
that $\tau_1(b)< \tau_2(b)$ and set $\tau_1(b)1=a$. Then
$\tau_1(a)=\tau_1(b)$ and $\tau_2(a)< \tau_2(b)$, and so for every bounded trace $\tau$,
$\tau(a)\leq \tau(b)$. This also holds for the unique unbounded trace with
value $\infty$ on all points except 0. But we cannot have $a\preceq b$, because
in that case $a+c\sim b$ for some $c\in A^+$, and so $\tau_1(c)=0$
and $\tau_2(c)>0$, which is impossible by the simplicity of $A$.

Let us consider the following topology on $T(A)$: 
the net of lower semicontinuous traces $(\tau_i)$ converges to $\tau$ if
\[
\limsup \tau_i((a-\epsilon)_+)\leq \tau(a)\leq \liminf \tau_i(a)
\]
for all $a\in A^+$ and $\epsilon>0$.
It was shown in \cite{elliott-robert-santiago} that  $T(A)$ is a
compact Hausdorff space and that the functor $T(\cdot)$ is a continuous contravariant functor 
from the category of C*-algebras to the category
of topological spaces.

For every positive element $a\in A^+$ the function $\widetilde a\colon T(A)\to
[0,\infty]$, defined by $\widetilde a(\tau)=\tau(a)$, is linear and lower semicontinuous. Consider the following 
question: Which linear and lower semicontinuous functions on $T(A)$ arise from positive
elements of $A$ in the form $\widetilde a$? In view of Theorem \ref{maincomparison},
this question is asking for
a description of the ordered semigroup $A_{CP}$ as a semigroup of functions on $T(A)$.
Theorem \ref{surjective} below gives an answer to this question assuming that $A$ is stable.

Let us recall the  definition of the non-cancellative cone $S(T(A))$ given in 
\cite{elliott-robert-santiago}. 
Let $LSC(T(A))$ denote the linear, lower semicontinuous, 
functions on $T(A)$ with values in $[0,\infty]$.
Then $S(T(A))$ is composed of the functions $f\in LSC(T(A))$
that satisfy the following condition: there is an increasing sequence $(h_n)$, $h_n\in LSC(T(A))$, such that
$f=\sup h_n$ and $h_n$ is continuous on each point where $h_{n+1}$ is finite.

For every $a\in A^+$ we have $\widetilde a\in S(T(A))$ (see \cite[Proposition 5.1]{elliott-robert-santiago}). 
Moreover, by \cite[Theorem 5.9]{elliott-robert-santiago}, $S(T(A))$ is the set of
functions $f\in LSC(T(A))$ for which there is
an increasing sequence of functions $\widetilde a_n$ coming from positive
elements $a_n\in A^+$ and such that $f=\sup \widetilde a_n$.
We can now prove the following theorem.

\begin{theorem}\label{surjective}
Let $A$ be a stable C*-algebra. Then $S(T(A))$ is the set of
functions of the form $\widetilde a$ for some $a\in A^+$.
\end{theorem}
\begin{proof}
Let $f$ be in $S(T(A))$ and let $(\widetilde{a_n})$, with $a_n\in A^+$, be an increasing
sequence of functions with supremum $f$. By Theorem \ref{maincomparison} (i),
we have that $a_n\precp a_{n+1}$ for all $n$. We may replace the positive elements
$a_n$ by $a_n'=(a_n-\epsilon_n)_+$, with $\epsilon_n>0$ small enough,  
so that $a_n'\preceq a_{n+1}'$ and $f=\sup \widetilde{a_n'}$.
Let $c_n\in A^+$ be such that $a_n'+c_n\sim a_{n+1}'$. Using the stability of $A$
we may find mutually orthogonal elements $c_n''$ such that $c_n''\sim c_n'$. Furthermore,
again by the stability of $A$, we may assume that the $c_n''$s have sufficiently small
norm such that the series $a=\sum_{n=1}^\infty c_n''$ converges. We then have $f=\widetilde{a}$. 
\end{proof}

\emph{Question 1.} Is $S(T(A)) = LSC(T(A))$?

\emph{Question 2.} In \cite{elliott-robert-santiago} a cone dual to the 2-quasitraces of $A$ is defined which
is analogous to the cone $S(T(A))$. Is the analog of Theorem \ref{surjective} for 2-quasitraces true?

\textbf{Acknowledgments.} I am grateful to George A. Elliott for pointing out that
Theorem \ref{surjective} was in all likelihood true.


\begin{bibdiv}
\begin{biblist}

  \bib{cuntz-pedersen}{article}{
   author={Cuntz, J.},
   author={Pedersen, G. K.},
   title={Equivalence and traces on C*-algebras},
   journal={J. Funct. Anal.},
   volume={33},
   date={1979},
   number={2},
   pages={135--164},
   issn={0022-1236},
}

  \bib{elliott-robert-santiago}{article}{
  author={Elliott, G. A.},
  author={Robert, L.},
  author={Santiago, L.},
  title={On the cone of lower semicontinuous traces on a C*-algebra},
  journal={arXiv:0805.3122v1 [math. OA]},
  date={2008},
}


\bib{kirchberg-rordam}{article}{
   author={Kirchberg, E.},
   author={R{\o}rdam, M.},
   title={Infinite non-simple C*-algebras: absorbing the Cuntz
   algebras $\scr O\sb \infty$},
   journal={Adv. Math.},
   volume={167},
   date={2002},
   number={2},
   pages={195--264},
   issn={0001-8708},
}

\bib{rordam}{article}{
   author={R{\o}rdam, M.},
   title={The stable and the real rank of $\scr Z$-absorbing C*-algebras},
   journal={Internat. J. Math.},
   volume={15},
   date={2004},
   number={10},
   pages={1065--1084},
   issn={0129-167X},
}

\end{biblist}
\end{bibdiv}
\end{document}